\documentclass[12pt]{amsart}%
\usepackage{latexsym}
\usepackage{amsthm}
\usepackage{amsmath}
\usepackage{amsfonts}
\usepackage{amssymb}
\usepackage[dvips]{graphicx}
\input xy
\xyoption{all}

\newtheorem{theo}{Theorem}[section]

\newtheorem{propo}[theo]{Proposition}
\newtheorem{defi}[theo]{Definition}
\newtheorem{coro}[theo]{Corollary}
\newtheorem{rem}[theo]{Remark}

\newtheorem*{rema}{Remark}
\newtheorem{examples}[theo]{Examples}

\newcommand\Ho{\operatorname{\it Ho}}
\newcommand\colim{\operatorname{\it colim}}
\newcommand\Ind{\operatorname{\it Ind}}

\newcommand\op{\operatorname{op}}

\newcommand\id{\operatorname{id}}
\newcommand\Id{\text{\it Id\,}}

\newcommand\hoco{\operatorname{\text{\it{hocolim}}}}
\newcommand\wco{\operatorname{\text{\it{wcolim}}}}

\newcommand\cc{\mathcal {C}}
\newcommand\cw{\mathcal {W}}
\newcommand\ck{\mathcal {K}}

\newcommand\cd{\mathcal {D}}
\newcommand\ca{\mathcal {A}}
\newcommand\cb{\mathcal {B}}
\newcommand\cx{\mathcal {X}}

\newcommand\cf{\mathcal {F}}

\newcommand\cl{\mathcal {L}}
\newcommand\ci{\mathcal {I}}
\newcommand\cj{\mathcal {J}}
\newcommand\cg{\mathcal {G}}

\newcommand\Set{\mathbf{Set}}

\newcommand\SSet{\mathbf{SSet}}
\newcommand\Ch{\mathbf{Ch}}

\newcommand\Ab{\mathbf{Ab}}

\newcommand\BD{\boldsymbol{\Delta}}

\newcommand\Spp{\mathbf{Sp}}

\font\aa=cmsy10 scaled 1300

\date{November 27, 2009}

\begin{document}
\title[Generalized Brown representability] 
{Generalized Brown representability in homotopy categories}
\author[Ji\v{r}\'{\i} Rosick\'{y}]{Ji\v{r}\'{\i} Rosick\'{y}$^*$}
\thanks{Supported by the Ministry of
Education of the Czech Republic under the project MSM 0021622409.}

\begin{abstract}
We show that the homotopy category of a combinatorial stable model category $\ck$ is well generated. It means that each object $K$
of $\Ho(\ck)$ is an iterated weak colimit of $\lambda$-compact objects for some cardinal $\lambda$. A natural question is whether
each $K$ is a weak colimit of $\lambda$-compact objects. We show that this is related to (generalized) Brown representability
of $\Ho(\mathcal K)$.
\end{abstract}

\maketitle
\section{Introduction}\label{sec0}
Combinatorial model categories were introduced by J. H. Smith as model categories which are locally presentable 
and cofibrantly generated. The latter means that both cofibrations and trivial cofibrations are cofibrantly generated 
by a set of morphisms. Most of important model categories are at least Quillen equivalent to a combinatorial one.
A natural question is to find properties of homotopy categories of combinatorial model categories. M. Hovey \cite{H}, 7.3.1
showed that $\Ho(\ck)$ has a set of weak generators for each cofibrantly generated pointed model category $\ck$. We will
prove that $\Ho(\ck)$ is well generated whenever $\ck$ is a combinatorial pointed model category. It means the existence
of a cardinal $\lambda$ such that $\Ho(\ck)$ has a set of weak $\lambda$-compact generators. This concept was introduced 
by A. Neeman \cite{N} for triangulated categories but it makes sense for homotopy categories of pointed model categories
as well. Our result generalizes that of \cite{H}, 7.4.3 proved for finitely generated pointed model categories.

If $\ck$ is a stable model category then the smallest localizing subcategory containing a set of weak generators $\cg$
is $\Ho(\ck)$ itself. It means that each object of $\Ho(\ck)$ is an iterated weak colimit of objects from $\cg$. We can ask
whether a set $\ca$ of objects of $\Ho(\ck)$ can be found such that each object of $\Ho(\ck)$ is a weak colimit of objects
from $\ca$. This question is related to (generalized) Brown representability of $\Ho(\ck)$. Consider the canonical functor
$$
E_\ca:\Ho(\ck)\to\Set^{\ca^{\op}}
$$
sending an object $K$ to the restriction
$$
E_\ca K =\hom (-, K)\big/\ca^{op}
$$
of its hom-functor $\hom(-, K): \Ho(\ck)^{\op}\to \Set$ to $\ca^{\op}$. Since $\Set^{\ca^{\op}}$ is a free completion
of $\ca$ under colimits (see \cite{AR}, 1.45), each object of $\Ho(\ck)$ is a weak colimit of objects from $\ca$ if and
only if $E_\ca$ is full. Let $\Ho(\ck)$ be compactly generated and $\ca$ consists of $\aleph_0$-compact objects. Then $E_\ca$ is
full if and only if $\Ho(\ck)$ satisfies [BRM] (see \cite{CKN}), i.e., iff it is Brown representable (for homology)
on morphisms. So, our question is whether $\Ho(\ck)$ satisfies [BRM$_\lambda$] for some cardinal $\lambda$. 
A. Beligiannis \cite{Be}, 11.8 showed that [BRM] implies [BRO], i.e., that every exact functor $\ca^{\op}\to\Ab$ ($\ca$
still consists of $\aleph_0$-compact objects) is in the image of $E_\ca$. Since these exact functors form the free completion
$\Ind(\ca)$ of $\ca$ under filtered colimits, [BRM] is equivalent with
$$
E_\ca:\Ho(\ck)\to\Ind(\ca)
$$
being full and surjective on objects. More precisely, one should say essentially surjective in the sense that each object 
from  $\Ind(\ca)$ is isomorphic to $E_\ca K$ for some $K$. 

Given a combinatorial stable model category $\ck$ such $\Ho(\ck)$ is well $\lambda$-generated for a regular cardinal $\lambda$, let $\ca$ denote 
the full subcategory consisting of $\lambda$-compact objects. Then the image of $E_\ca$ is contained in the free completion $\Ind_\lambda(\ca)$ 
of $\ca$ under $\lambda$-filtered colimits. Our generalized Brown representability thus means the question whether
$$
E_\ca:\Ho(\ck)\to\Ind_\lambda(\ca)
$$
is full and essentially surjective on objects. Previous versions of this paper claimed that for each combinatorial stable model category $\ck$
there is a regular cardinal $\lambda$ such that this is true. Unfortunately, the proofs contain a gap and the author is grateful to R. Jardine
and F. Muro for pointing this up. Let us add the one cannot expect $E_\ca$ being also faithfull. Then $\Ho(\ck)$ would be equivalent
to $\Ind_\lambda(\ca)$, i.e., it would be accessible. Even in the compactly generated case, i.e., for $\lambda=\aleph_0$, there 
are, except trivial situations, phantoms, i.e., non-zero morphisms $f$ in $\Ho(\ck)$ with $E_\ca(f)=0$.  

\section{Basic concepts}\label{sec1}

A \textit{model structure} on a category $\ck$ will be understood
in the sense of Hovey \cite{H}, i.e., as consisting of three
classes of morphisms called weak equivalences, cofibrations and
fibrations which satisfy the usual properties of Quillen \cite{Q}
and, moreover, both (cofibration, trivial fibrations) and
(trivial cofibrations, fibration) factorizations are functorial.
Recall that trivial (co)fibrations are those (co)fibrations which
are in the same time weak equivalences. The (cofibration, trivial
fibration) factorization is \textit{functorial} if there is a functor
$F:\ck^{\to}\to\ck$ and natural transformations $\alpha:\mbox{dom}
\to F$ and $\beta:F\to\mbox{cod}$ such that $f=\beta_f\alpha_f$
is the (cofibration, trivial fibration) factorization of $f$.
Here $\ck^{\to}$ denotes the category of morphisms in $\ck$
and $\mbox{dom}:\ck^{\to}\to\ck$ ($\mbox{cod}:\ck^{\to}\to\ck$)
asssign to each morphism its (co)domain. The same for
(trivial cofibration, fibration) factorization (see \cite{RT}) .

A \textit{model category}
is a complete and cocomplete category together with a model
structure. In a model category $\ck$, the classes of weak equivalences,
cofibrations and fibrations will be denoted by $\cw$, $\cc$ and
$\cf$, resp. Then $\cc_0 = \cc \cap \cw$ and $\cf_0 =\cf\cap \cw$
denote trivial cofibrations and trivial fibrations, resp. We have
$$
\cf_0 =\cc^{\square}\,, \quad \cf=\cc_0^{\square}\,, \quad \cc
={}^\square \cf_0\quad\mbox{and}\quad \cc_0={}^{\square}\cf
$$
where $\cc^{\square}$ denotes the class of all morphisms having
the right lifting property w.r.t.\ each morphism from $\cc$ and
${}^\square \cf$ denotes the class of all morphisms having the left
lifting property w.r.t.\ each morphism of $\cf$. $\ck$ is called
\textit{cofibrantly generated} if there are sets of morphisms
$\ci$ and $\cj$ such that $\cf_0 = \ci^{\square}$ and $\cf =
\cj^{\square}$.  If $\ck$ is locally presentable then $\cc$ is
the closure of $\ci$ under pushouts, transfinite compositions and
retracts in comma-categories $K\downarrow \ck$ and, analogously,
$\cc_0$ is this closure of $\cj$.

An object $K$ of a model category $\ck$ is called
\textit{cofibrant} if the unique morphism $0\to K$ from an
initial object is a cofibration and $K$ is called \textit{fibrant}
if the unique morphism $K\to 1$ to a terminal object is a
fibration. Let $\ck_c$, $\ck_f$ or $\ck_{cf}$ denote the full
subcategories of $\ck$ consisting of objects which are
cofibrant, fibrant or both cofibrant and fibrant resp.
We get the \textit{cofibrant replacement functor} $R_c:\ck\to\ck$
and the \textit{fibrant replacement functor} $R_f:\ck\to\ck$.
We will denote by $R=R_fR_c$ their composition and call it the
\textit{replacement functor}. The codomain restriction of 
the replacement functors are $R_c:\ck\to\ck_c$, $R_f:\ck\to\ck_f$
and $R:\ck\to\ck_{cf}$.

Let $\ck$ be a model category and $K$ an object of $\ck$. Recall
that a \textit{cylinder object} $C(K)$ for $K$ is given by a
(cofibration, weak equivalence) factorization
$$
\nabla : K\amalg K \xrightarrow{\ \gamma_K\ } C(K)\xrightarrow{\
\sigma_K\ } K
$$
of the codiagonal $\nabla$. Morphisms $f,g:K\to L$ are
\textit{left homotopic} if there is a morphism $h: C(K)\to L$
with
$$
f= h\gamma_{1K} \quad \mbox{and}\quad g=h\gamma_{2K}
$$
where $\gamma_{1K}= \gamma_K i_1$ and $\gamma_{2K}=\gamma_K
i_2$ with $i_1,i_2 :K\to K\amalg K$ being the coproduct
injections. In fact, cylinder objects form a part of the \textit{cylinder 
functor} $C:\ck\to\ck$ and $\gamma_1,\gamma_2:\Id\to C$ are natural transformations.

On $\ck_{cf}$, left homotopy $\sim$ is an equivalence
relation compatible with compositions, it does not
depend on a choice of a cylinder object and we get the quotient
$$
Q :\ck_{cf} \to \ck_{cf}/\sim \,.
$$
The composition
$$
P:\ck\xrightarrow{\ R\ } \ck_{cf}
\xrightarrow{\ Q\ } \ck_{cf}/\sim
$$
is, up to equivalence, the projection of $\ck$ to the homotopy
category $\Ho(\ck)= \ck[\cw^{-1}]$ (see \cite{H}). In what
follows, we will often identify $\ck_{cf}/\sim$ with $\Ho(\ck)$.

A category $\ck$ is called $\lambda$-\textit{accessible}, where $\lambda$ is
a regular cardinal, provided that

(1) $\ck$ has $\lambda$-filtered colimits,

(2) $\ck$ has a set $\ca$ of $\lambda$-presentable objects such that
every object\newline
\indent\hskip 18pt of $\ck$ is a $\lambda$-filtered colimit of objects from
$\ca$.

\noindent
Here, an object $K$ of a category $\ck$ is called
$\lambda$-\textit{presentable} if its hom-functor $\hom(K,-):\ck\to\Set$ preserves
$\lambda$-filtered colimits; $\Set$ is the category of sets. A category is
called \textit{accessible} if it is $\lambda$-accessible for some regular cardinal
$\lambda$. The theory of accessible categories was created in \cite{MP} and for its
presentation one can consult \cite{AR}. We will need to know that $\lambda$-accessible
categories are precisely categories $\Ind_\lambda(\ca)$ where $\ca$ is a small category.
If idempotents split in $\ca$ then $\ca$ precisely consists of $\lambda$-presentable
objects in $\Ind(\ca)$. In what follows, we will denote by $\ck_\lambda$ the full
subcategory of $\ck$ consisting of $\lambda$-presentable objects.

A \textit{locally $\lambda$-presentable category} is defined as a cocomplete
$\lambda$-acce\-ssible category and it is always complete. Locally
$\lambda$-presentable ca\-te\-go\-ries are precisely categories
$\Ind_\lambda(\ca)$ where
the category $\ca$ has $\lambda$-{\textit{small}} co\-li\-mits, i.e.,
colimits of
diagrams $D:\cd\to\ca$ where $\cd$ has less then $\lambda$ morphisms.
In general,
the category $\Ind_\lambda(\ca)$ can be shown to be the full subcategory of
the functor category
$\Set^{\ca^{op}}$ consisting of $\lambda$-filtered colimits $H$ of hom-functors
$\hom(A,-)$ with $A$  in $\ca$. In the case that $\ca$ has $\lambda$-small colimits
this is equivalent to the fact that  $H:\ca^{op}\to\Set$ preserves
$\lambda$-small limits.
More generally, if $\ca$ has weak $\lambda$-small co\-li\-mits then
$\Ind_\lambda(\ca)$
precisely consists of left $\lambda$-covering functors (see \cite{Hu} 3.2).
Let us recall that a weak colimit of a diagram $D:\cd\to\ca$ is a cocone
from $D$
such that any other cocone from $D$ factorizes through it but
not necessarily uniquely.
If $\cx$ is a category with weak $\lambda$-small limits then a functor $H:\cx\to\Set$
is \textit{left $\lambda$-covering} if, for each $\lambda$-small diagram
$D:\cd\to\cx$ and its weak limit $X$, the canonical mapping $H(X)\to\lim HD$
is surjective (see \cite{CV} for $\lambda=\omega$). A left $\lambda$-covering
functor preserves all $\lambda$-small limits which exist in $\cx$. Moreover,
a functor $H:\cx\to\Set$ is left $\lambda$-covering iff it is weakly
$\lambda$-continuous, i.e., iff it preserves weak $\lambda$-small limits.
This immediately follows from \cite{CV}, Proposition 20 and the fact that
surjective mappings in $\Set$ split. A functor $H$ is called
\textit{weakly continuous} if it preserves weak limits. Hence a weakly
continuous functor $H:\cx\to\Set$ preserves all existing limits.

A functor $F:\ck\to\cl$ is called $\lambda$-\textit{accessible} if $\ck$ and $\cl$ are
$\lambda$-acce\-ssible categories and $F$ preserves $\lambda$-filtered colimits. An
important subclass of $\lambda$-accessible functors are those functors which also
preserve $\lambda$-presentable objects. In the case that idempotents split in $\cb$, 
those functors are precisely functors $\Ind_\lambda(G)$ where $G:\ca\to\cb$ is a functor. 
The uniformization theorem of Makkai and Par\' e says that for each $\lambda$-acce\-ssible functor 
$F$ there are arbitrarily large regular cardinals $\mu$ such that $F$ is $\mu$-accessible and 
preserves $\mu$-presentable objects (see \cite{AR} 2.19). In fact, one can take 
$\lambda\triangleleft\mu$ where $\triangleleft$ is the set theoretical relation between regular
cardinals corresponding to the fact that every $\lambda$-accessible category is $\mu$-accessible 
(in contrast to \cite{AR} and \cite{MP}, we accept $\lambda\triangleleft\lambda$). For every 
$\lambda$ there are arbitrarily large regular cardinals $\mu$ such that
$\lambda\triangleleft\mu$. For instance, $\omega\triangleleft\mu$ for every regular cardinal
$\mu$.  

\section{Combinatorial model categories}\label{sec2}

We will follow J.~H.~Smith and call a model category $\ck$
$\lambda$-\textit{combina\-to\-rial} if $\ck$ is locally $\lambda$-presentable
and both cofibrations and trivial cofibrations are cofibrantly generated
by sets $\ci$ and $\cj$ resp. of morphisms having $\lambda$-presentable domains
and codomains. Then both tri\-vi\-al fibrations and fibrations are closed in $\ck^\to$
under $\lambda$-filtered colimits. $\ck$ will be called \textit{combinatorial}
if it is $\lambda$-combinatorial for some regular cardinal $\lambda$. Clerly,
if $\lambda<\mu$ are regular cardinals and $\ck$ is $\lambda$-combinatorial then
$\ck$ is $\mu$-combinatorial. 

The following result is due to J.~H.~Smith  and is presented in \cite{D$_1$}, 7.1. 
We just add a little bit more detail to the proof.

\begin{propo}[Smith]\label{prop2.1}
Let $\ck$ be a $\lambda$-combinatorial model category. Then the functors
$\ck^{\to}\to\ck$ giving (cofibration, trivial fibration) and
(trivial cofibration, fibration) factorizations are $\lambda$-accessible.
\end{propo}

\begin{proof}
We know that $\ck$ is locally $\lambda$-presentable and domains and co\-do\-mains of morphisms from the generating 
set $\ci$ of cofibrations are $\lambda$-presentable. For every morphism $f:A\to B$ form a colimit $F_0f$ of the diagram
$$
\xymatrix{
A \\
&&\\
X\ar[uu]^u  \ar@{.}[ur] \ar [rr]_h && Y
}
$$
consisting of all spans $(u,h)$ with $h:X\to Y$ in $\ci$ such that
there is $v:Y\to B$ with $vh=fu$. Let $\alpha_{0f}:A\to F_0f$ denote
the component of the colimit cocone (the other components are $Y\to F_0f$
and they make all squares
$$
\xymatrix@C=3pc@R=3pc{
A \ar[r]^{\alpha_{0f}} & F_0f \\
X \ar [u]^{u} \ar [r]_{h} &
Y \ar[u]_{}
}
$$
to commute). Let $\beta_{0f}:F_0f\to B$ be the morphism induced by $f$ and $v$'s. Then 
$F_0:\ck^{\to}\to\ck$ is clearly $\lambda$-accessible. Let $F_if,\alpha_{if}$ and
$\beta_{if}$, $i\leq\lambda$, be given by the following transfinite
induction: $F_{i+1}f=F_0\beta_{if}$, $\alpha_{i+1,f}=\alpha_{0,\beta_{if}}
\alpha_{if}$, $\beta_{i+1,f}=\beta_{0,\beta_{if}}$ and the limit step
is given by taking colimits. Then all functors $F_i:\ck^{\to}\to\ck$,
$i\leq\lambda$ are $\lambda$-accessible and $F_\lambda$ yields the desired
(cofibration, trivial fibration) factorization.

The proof for (trivial cofibration, fibration) factorizations is analogous.
\end{proof}

\begin{rem}\label{re2.2}
{\em
(1) \ref{prop2.1} implies that, in a $\lambda$-combinatorial model category $\ck$, weak equivalences are closed 
under $\lambda$-filtered colimits in $\ck^\to$ (see \cite{D$_1$}, 7.5).

(2) Following the uniformization theorem (\cite{AR} Remark 2.19), there is a regular cardinal
$\mu$ such that the functors from \ref{prop2.1} are $\mu$-accessible and preserve
$\mu$-presentable objects. This means that the factorizations $A\to C\to B$ of a morphism
$A\to B$ have $C$ $\mu$-presentable whenever $A$ and $B$ are $\mu$-presentable. This
point is also well explained in \cite{D$_1$}, 7.2.
}
\end{rem}

\begin{defi}\label{defi2.3}
{\em
A $\lambda$-combinatorial model category $\ck$ will be called \textit{strongly $\lambda$-combinatorial}
if the functor $F:\ck^{\to}\to\ck$ giving the (cofibration, trivial fibration) factorization preserves 
$\lambda$-presentable objects.
}
\end{defi}

\begin{rem}\label{re2.4}
{\em
(1) Following \ref{re2.2} (2), every combinatorial model category is strongly $\mu$-combinatorial for some
regular cardinal $\mu$. 

(2) Following \cite{AR}, 2.20, if $\ck$ is strongly $\lambda$-combinatorial and $\lambda\triangleleft\mu$ 
then $\ck$ is strongly $\mu$-combinatorial.

(3) In a strongly $\lambda$-combinatorial model category $\ck$, both the co\-fib\-rant replacement functor 
$R_c:\ck\to\ck$ and the cylinder functor $C:\ck\to\ck$ preserve $\lambda$-filtered colimits and $\lambda$-presentable 
objects.
}
\end{rem}

Combinatorial model categories form a very broad class. We will give some important examples.  

\begin{examples}\label{exa2.5}
{\em
(i) The model category $\SSet$ of simplicial sets is strongly $\omega_1$-combinatorial. Let us add that the functor
$F:\ck^\to\to\ck$ giving (trivial cofibration, fibration) factorization preserves $\omega_1$-pre\-sen\-tab\-le objects
too. This observation can be found in \cite{I}, Section 5. 

The same is true for the model category $\SSet_\ast$ of pointed simplicial sets.

(ii) The category $\Ch(R)$ of chain complexes of modules over a ring $R$ is an $\omega$-combinatorial model 
category (see \cite{H}, 2.3.11); fibrations are dimensionwise surjections and weak equivalences are homology 
isomorphisms. We will show that this model category is strongly $\omega$-com\-bi\-na\-to\-rial provided that $R$ 
is a noetherian ring of a finite projective dimension.

Finitely presentable objects in $\Ch(R)$ are precisely bounded complexes of finitely presentable modules. 
Each such a complex is clearly finitely presentable. On the other hand, consider a chain complex $(A,d)$;
it means that $(d_n:A_n\to A_{n-1})$ for each integer $n$. For each integer $k\geq 0$, let $(A^k,d^k)$ be
the following chain complex: $A^k_n=0$ for $n>k$ and $n<-k-1$, $A^k_n=A_n$ for $-k\leq n\leq k$, 
$A^k_{-k-1}=A_{-k}$, $d^k_n=0$ for $n\leq -k$ and $n>k$, $d^k_n=d_n$ for $-k<n\leq k$ and $d^k_{-k}=\id_{A_{-k}}$. 
Then $(A,d)$ is a colimit of the chain $(A_k,d^k)$ with colimit components $f$ where $f_n=0$ for $n<-k-1$ and $n>k$, 
$f_n=\id_{A_n}$ for $-k\leq n\leq k$ and $f_{-k-1}=d_{-k}$. Thus each finitely presentable complex is bounded
and evidently consists of finitely presentable modules.

Following \cite{CH}, 2.9, $\Ch(R)$ is strongly $\omega$-combinatorial if and only if finitely presentable
complexes have finitely presentable cofibrant replacements.
 
(iii) The category $\Spp$ of spectra with the strict model category structure (in the sense of \cite{BF}) 
is $\omega$-combinatorial (see \cite{S} A.3). We will show that it is strongly $\omega_1$-combinatorial.
 
Let us recall that a spectrum $X$ is a sequence $(X_n)_{n=0}^\infty$
of pointed simplicial sets equipped with morphisms $\sigma_n^X
:\Sigma X_n\to X_{n+1}$ where $\Sigma$ is the suspension functor.
This means that $\Sigma X_n = S^1\wedge X_n$ where $S^1\wedge -$
is the smash product functor, i.e., a left adjoint to
$$
-{}^{S^1}=\hom(S^1, -):\SSet_\ast\to\SSet_\ast.
$$

The strict model structure on $\Spp$ has level equivalences as weak equivalences
and level fibrations as fibrations. This means that
$f:X\to Y$ is a weak equivalence (fibration) iff all
$f_n:X_n\to Y_n$ are weak equivalences (fibrations) in
$\SSet_\ast$. A morphism $f:X\to Y$ is a (trivial) cofibration iff
$f_0:X_0\to Y_0$ is a (trivial) cofibration and all induced
morphisms $t_n : Z_n\to Y_n$, $n\geq 1$, from pushouts are
(trivial) cofibrations
$$
\xymatrix@C=3pc{
\Sigma X_{n-1}\ar [r]^{\sigma_{n-1}^X}
              \ar [d]_{\Sigma f_{n-1}}&
X_n \ar [d] \ar [ddr]^{f_n}&\\
\Sigma Y_{n-1} \ar [r]
               \ar [drr]_{\sigma_{n-1}^Y}&
Z_n \ar [dr]^{t_n}&\\
&& Y_n
}
$$
(see \cite{BF}, \cite{HSS} or \cite{H$_1$}). Then a (cofibration,
trivial fibration) factorization $X\xrightarrow{\ g\ }Z
\xrightarrow{\ h\ } Y$ of a morphism $f:X\to Y$ is  made as
follows.

One starts with a (cofibration, trivial fibration) factorization
$$
f_0 : X_0 \xrightarrow{\ g_0\ } Z_0 \xrightarrow{\ h_0\ } Y_0
$$
in $\SSet_\ast$. Then one takes a (cofibration, trivial fibration)
factorization
$$
t:Z_1^\prime  \xrightarrow{\ u\ } Z_1 \xrightarrow{\  h_1\ } Y_1
$$
of the induced morphism from a pushout
$$
\xymatrix@C=3pc{
\Sigma X_{0}\ar [r]^{\sigma_{0}^X}
              \ar [d]_{\Sigma g_{0}}&
X_1 \ar [d]_{q}
    \ar [ddr]^{f_1}&\\
\Sigma Z_{0} \ar [r]^{p}
               \ar [drr]_{\sigma_{0}^Y\cdot \Sigma h_0}&
Z^\prime_1 \ar [dr]_<<{t}&\\
&& Y_1
}
$$
and puts $\sigma_1^Z=up$ and $g_1=uq$. This yields
$$
f_1:X_1 \xrightarrow{\ g_1\ } Z_1 \xrightarrow{\ h_1\ } Y_1
$$
and one continues the procedure. Analogously, one constructs a (trivial cofibration, fibration) factorization. Since
a spectrum $X$ is $\omega_1$-presentable iff all $X_n$, $n\geq 0$ are $\omega_1$-presentable in $\SSet_\ast$, it is 
easy to see that the strict model structure on $\Spp$ is strongly $\omega_1$-combinatorial. Moreover, the functor
$F:\ck^\to\to\ck$ giving (trivial cofibration, fibration) factorization preserves $\omega_1$-presentable objects too. 

(iv) The model category $\Spp$ of spectra with the stable
Bousfield-Fried\-lan\-der model category structure (see \cite{BF}) is
$\omega$-combinatorial (see \cite{S} A.3). The stable model structure is defined 
as a Bousfield localization of the strict model structure, i.e., by adding a set of new
weak equi\-va\-len\-ces. Cofibrations and trivial fibrations remain
unchanged, which means that the stable model category of spectra is strongly $\omega_1$-combinatorial.
}
\end{examples} 

There is well known that the homotopy category of any model category $\ck$ has products, coproducts, weak limits
and weak colimits. We will recall their constructions.

\begin{rem}\label{rem2.6}
{\em
(i) Let $K_i$, $i\in I$ be a set of objects of $\ck$. Without any loss
of generality, we may assume that they are in $\ck_{cf}$. Then their product
in $\ck$
$$
p_i : K \to K_i
$$
is fibrant and let
$$
q_K : R_cK \to K
$$
be its cofibrant replacement. Then $R_cK \in  \ck_{cf}$ and
$$
Q(p_i q_K) : QR_cK \to QK_i
$$
is a product in $\Ho(\ck)$. Recall that $Q:\ck_{cf}\to\Ho(\ck)=\ck_{cf}/\sim$
is the quotient functor.

In fact, consider morphisms
$$
Q f_i : QL \to QK_i\,,\quad\quad i\in I
$$
in  $\ck_{cf}/\sim$. Let $f:L\to K$ be the induced morphism and
$g:L\to R_cK$ be given by the lifting property:
$$
\xymatrix@C=3pc@R=3pc{
0 \ar [r] \ar[d] & R_cK \ar[d]^{q_K}
\\
L\ar @{-->}[ur]_g  \ar[r]_f & K
}
$$
We have $Q(p_i q_K g)= Qf_i$ for each $i\in I$. The
unicity of $g$ follows from the facts that $Qq_K$ is an isomorphism
and that left homotopies $h_i$ from $p_i f$ to $p_i f'$, $i\in
I$, lift to the left homotopy from $f$ to $f'$.

Since $\ck^{\op}$ is a model category and
$$
\Ho (\ck^{\op})= (\Ho (\ck))^{\op}\,,
$$
$\Ho (\ck)$ has coproducts.

\vskip 1mm
(ii) In order to show that $\Ho(\ck)$ has weak colimits, it
suffices to prove that it has weak pushouts.
In fact, a weak coequalizer
$$
\xymatrix@1{
A\ \ar@<0.6ex>[r]^{f}\
 \ar@<-0.6ex>[r]_g\ &\ B\ \ar [r]^h\ &\ D
}
$$
is given by a weak pushout
$$
\xymatrix@C=3pc@R=3pc{
B \ar [r]^h  & D\\
A\amalg B \ar [u]^{(f, \id_B)} \ar [r]_{(g, \id_B)} &
B\ar [u]_h
}
$$
and weak colimits are constructed using coproducts and weak
co\-equa\-li\-zers in the same way as colimits are constructed  by
coproducts and coequalizers. This means that, given a diagram $D:\cd\to\Ho(\ck)$, its weak colimit $K$
is a weak coequalizer of $f$ and $g$ constructed as follows
$$
\xymatrix@C=3pc@R=3pc{
Dd \ar[d]_{u_e}  \ar[rd]^{v_d}& &\\
\coprod\limits_{e:d\to d'} Dd \ar@<1.2ex>[r]^f \ar@<-0.2ex>[r]_g &
\coprod\limits_d Dd \ar@<0.8ex>[r]^h&\underset{\ }{K}\\
Dd \ar[u]^{u_e} \ar[r]_{De} &Dd' \ar[u]_{v_{d'}}&
}
$$
where $u_e$ and $v_d$ are coproduct injections. The weak colimit cocone $\delta_d:Dd\to K$ is given by
$$
\delta_d=hv_d
$$
for each $d$ in $\cd$. We emphasize that the coproduct on the left is over all morphisms of $\cd$.

Let
$$
\xymatrix@C=3pc@R=3pc{
B &\\
P\ar[u]^{f} \ar[r]_{g} & D
}
$$
be a diagram in $\ck$. Consider a pushout
$$
\xymatrix@C=3pc@R=3pc{
B_1 \ar[r]^{\overline{g}} & E \\
A \ar [u]^{f_1} \ar [r]_{g_1} &
D_1 \ar[u]_{\overline{f}}
}
$$
in $\ck$ where $f=f_2f_1$ and $g=g_2g_1$ are (cofibration, trivial fibration)
factorizations. Following the homotopy extension property of cofibrations (see \cite{Hi}, 7.3.12),
$$
\xymatrix@C=3pc@R=3pc{
PB_1 \ar[r]^{P\overline{g}} & PE \\
PA \ar [u]^{Pf_1} \ar [r]_{Pg_1} &
PD \ar[u]_{P\overline{f}}
}
$$
is a weak pushout in $\Ho(\ck)$ which is called the \textit{homotopy pushout}
of the starting diagram. Recall that $P:\ck\to\Ho(\ck)$ is the canonical fun\-ctor.

Following \cite{C}, we will call the resulting weak colimits in $\Ho(\ck)$
\textit{standard}. By duality, $\Ho(\ck)$ has weak limits. Since our model categories are functorial,
the construction in $\ck$ giving standard weak colimits in $\Ho(\ck)$ is functorial in $\ck$.

(iii) Consider a diagram $D:\cd\to \ck$, its colimit $(\bar
\delta_d : Dd \to \overline K)$ and $(\delta_d : Dd \to K)$
such that $(P\delta_d : PDd \to PK)$ is a
standard weak colimit of $PD$. There is the comparison morphism
$p:K\to \overline K$  such that $P(k)\delta_d = P(\overline\delta_d)$ for
each $d\in \cd$. It suffices to find  this morphism for a 
pushout diagram
$$
\xymatrix@C=3pc@R=3pc{
B \ar[r]^{g'} & \overline E \\
A \ar [u]^{f} \ar [r]_{g} &
D \ar[u]_{f'}
}
$$
But it is given by $p\overline g=g'f_2$ and $p\overline f=f'g_2$; we use
the notation from (ii).

(iv) Another, and very important, colimit construction in model categories are
\textit{homotopy colimits} (see, e.g., \cite{BK}, \cite{DHKS}, \cite{Hi}). 
Both coproducts and homotopy pushouts described above are instances of this concept.
While weak colimits correspond to homotopy commutative diagrams, homotopy colimits
correspond to homotopy coherent ones. So, one cannot expect that homotopy colimits are weak
colimits. There is always a morphism $\wco D\to\hoco D$ from the standard weak colimit
to the homotopy colimit for each diagram $D:\cd\to\ck$. 

Following \ref{re2.2}(1), $\lambda$-filtered colimits are homotopy $\lambda$-filtered colimits 
in a $\lambda$-combinatorial model category.
}
\end{rem}

\section{Well generated homotopy categories}\label{sec4}

Given a small, full subcategory $\ca$ of a category $\ck$, the
\textit{canonical functor}
$$
E_\ca: \ck \to \Set^{\ca^{\op}}
$$
assigns to each object $K$ the restriction
$$
E_\ca K =\hom (-, K)\big/\ca^{op}
$$
of its hom-functor $\hom(-, K): \ck^{\op}\to \Set$ to $\ca^{\op}$
(see \cite{AR} 1.25). This functor is (a) $\ca$-\textit{full}
and (b) $\ca$-\textit{faithful} in the sense that

(a) for every $f:E_\ca A\to E_\ca K$ with $A$ in $\ca$ there is
$f':A\to K$\newline\indent{\hskip 15pt} such that $E_\ca f'=f$ and

(b) $E_\ca f = E_\ca g$ for $f,g:A\to K$ with $A$ in $\ca$ implies $f=g$.
\vskip 3 mm

Let $\ck$ be a locally $\lambda$-presentable model category and
denote by $\Ho(\ck_\lambda)$ the full subcategory $P(\ck_\lambda)$ of
$\Ho(\ck)$ consisting of $P$-images of $\lambda$-presentable
objects in $\ck$ in the canonical functor $P:\ck\to\Ho(\ck)$.
Let 
$$
E_\lambda:\Ho(\ck)\to \Set^{\Ho(\ck_\lambda)^{\op}}
$$ 
denote the canonical functor $E_{\Ho(\ck_\lambda)}$.

\begin{theo}\label{th3.1}
Let $\ck$ be a strongly $\lambda$-combinatorial model category. Then the composition 
$$
E_\lambda P:\ck\to\Set^{\Ho(\ck_\lambda)^{\op}}
$$ 
preserves $\lambda$-filtered colimits.
\end{theo}

\begin{proof}
Consider a $\lambda$-filtered diagram $D:\cd \to \ck$ and its colimit $(k_d : Dd \to K)$ in $\ck$. Consider 
$X\in \ck_\lambda$ and a morphism $f:PX\to PK$ in $\Ho(\ck)$. Let $u_X:R_cX\to X$ denote the cofibrant replacement
and $v_K:K\to R_fK$ the fibrant replacement. Following \cite{H}, 1.2.10(ii), there is $\overline{f}:R_cX\to R_fK$
such that $P(v_K)fPu_X=P\overline{f}$. Since $\ck$ is strongly $\lambda$-combinatorial, the object $R_cX$ is  
$\lambda$-presentable and $(R_fk_d : R_fDd \to R_fK)$ is a $\lambda$-filtered colimit. Thus $\overline{f}=R_f(k_d)g$
for some $g:R_cX\to R_fDd$ and $d$ in $\cd$. We have
$$
PR_f(k_d)P(g)(Pu_X)^{-1}=P(\overline{f})(Pu_X)^{-1}=P(v_K)f
$$
and thus
$$
f=(Pv_K)^{-1}PR_f(k_d)P(g)(Pu_X)^{-1}.
$$
Since
$$
v_Kk_d=R_f(k_d)v_{Dd}
$$
where $v_{Dd}:Dd\to R_fDd$ is the fibrant replacement of $Dd$, we have
$$
f=P(k_d)(Pv_{Dd})^{-1}P(g)(Pu_X)^{-1}.
$$
This proves that $f$ factorizes through some $Pk_d$. In order to verify that 
$E_\lambda Pk_d:E_\lambda PDd\to E_\lambda PK$ is a $\lambda$-filtered colimit, we have to show that
this factorization is essentially unique.

Assume that $f=P(k_d)g_1=P(k_d)g_2$ are two such factorizations, i.e., $g_1,g_2:PX\to PDd$. Again, using
\cite{H}, 1.2.10(ii), there are $\overline{g_i}:R_cX\to R_fDd$ such that $P\overline{g}_i=P(v_{Dd})g_iPu_X$
for $i=1,2$. Since $PR_f(k_d)P\overline{g}_1=PR_f(k_d)P\overline{g}_2$, the morphisms $R_f(k_d)\overline{g}_1$
and $R_f(k_d)\overline{g}_2$ are left homotopic (see \cite{H}, 1.2.10(ii) and 1.2.6). Thus there is a morphism
$h:CR_cX\to R_fK$ such that $R_f(k_d)\overline{g}_i=h\gamma_{iR_cX}$ for $i=1,2$. Since $CR_cX$ is 
$\lambda$-presentable, we can assume without any loss of generality that $h=R_f(k_d)\overline{h}$  
for $\overline{h}:CR_cX\to R_fDd$. Since
$$
R_f(k_d)\overline{h}\gamma_{1R_cX}=R_f(k_d)\overline{h}\gamma_{2R_cX}
$$
and $R_cX$ is $\lambda$-presentable, there is $e:d\to d'$ in $\cd$ such that
$$
R_fD(e)\overline{h}\gamma_{1R_cX}=R_fD(e)\overline{h}\gamma_{2R_cX}.
$$
Thus $R_fD(e)\overline{g}_1$ and $R_fD(e)\overline{g}_2$ are left homotopic. Therefore 
$$
PR_fD(e)P\overline{g}_1=PR_fD(e)P\overline{g}_2.
$$ 
Hence 
$$
PR_fD(e)P(v_{Dd})g_1Pu_X=PR_fD(e)P(v_{Dd})g_2Pu_X
$$
and thus 
$$
P(v_{Dd'})PD(e)g_1Pu_X=P(v_{Dd'})PD(e)g_2Pu_X.
$$
Consequently
$$
PD(e)g_1=PD(e)g_2,
$$
which is the desired essential unicity of our factorization.
\end{proof}

Let $P_\lambda : \ck_\lambda\to \Ho(\ck_\lambda)$ denote the
domain and codomain restriction of the canonical functor
$P:\ck\to\Ho(\ck)$. We get the induced functor
$$
\Ind_\lambda P_\lambda :\ck = \Ind_\lambda \ck_\lambda \to
\Ind_\lambda \Ho(\ck_\lambda)\,.
$$

\begin{coro}\label{cor3.2}
Let $\ck$ be a strongly $\lambda$-combinatorial model category. Then $E_\lambda P\cong \Ind_\lambda P_\lambda$.
\end{coro}

\begin{rema}
{\em
This means that $E_\lambda$ factorizes through the inclusion
$$
\Ind_\lambda \Ho(\ck_\lambda) \subseteq \Set^{\Ho(\ck_\lambda)^{\op}}
$$
and that the codomain restriction of $E_\lambda$, which we denote
$E_\lambda$ as well, makes the composition $E_\lambda P$
isomorphic to $\Ind_\lambda P_\lambda$.
}
\end{rema}

\begin{proof}
Since both $E_\lambda P$ and $\Ind_\lambda P_\lambda$ have the
same domain restriction on $\ck_\lambda$, the result follows from
\ref{th3.1}.
\end{proof}

\begin{coro}\label{cor3.3}
Let $\ck$ be a strongly $\lambda$-combinatorial model category. The the functor
$$
E_\lambda:\Ho(\ck_\lambda)\to\Ind_\lambda\Ho(\ck_\lambda)
$$
preserves coproducts.
\end{coro}
\begin{proof}
Following \ref{rem2.6} (i) and \ref{cor3.2}, it suffices to show that $\Ind_\lambda P_\lambda$ preserves coproducts.
Since each coproduct is a $\lambda$-filtered colimit of $\lambda$-small coproducts and $\Ind_\lambda P_\lambda$
preserves $\lambda$-filtered colimits, we have to prove that $\Ind_\lambda P_\lambda$ preserves $\lambda$-small 
coproducts. Let $\coprod\limits_{i\in I} K_i$ be such a coproduct, i.e., $card I<\lambda$. Each $K_i$ is 
a $\lambda$-filtered colimit $\colim D_i$ of $\lambda$-presentable objects. Let $D_i:\cd_i\to\ck_\lambda$ denote
the corresponding diagrams. Since $\coprod\limits_{i\in I}\colim D_i$ is isomorphic to a $\lambda$-filtered colimit 
of coproducts $\coprod\limits_{i\in I} D_id_i$ where $d_i\in\cd_i$, $\Ind_\lambda P_\lambda$ preserves $\lambda$-filtered 
colimits and $P_\lambda$ preserves $\lambda$-small coproducts, the result is proved.
\end{proof}

\begin{defi}\label{defi3.4}
{\em
Let $\ck$ be a category with coproducts and $\lambda$ a cardinal. An object $A$ of $\ck$ is called
$\lambda$-\textit{small} if for every morphism $f:A\to\coprod\limits_{i\in I}K_i$ there is a subset $J$ of $I$
of cardinality less than $\lambda$ such that $f$ factorizes as
$$
A\to\coprod\limits_{j\in J}K_j\to\coprod\limits_{i\in I}K_i
$$   
where the second morphism is the subcoproduct injection.
}
\end{defi}
 
\begin{rem}\label{re3.5}
{\em
$\aleph_0$-small objects are also called compact or abstractly finite. We use the terminology of A. Neeman \cite{N}
who found how compactness should be defined for uncountable cardinals. His definition was simplified by H. Krause
in \cite{K2}. They considered compactness in additive categories but the definition makes sense in general.

Consider classes $\mathcal S$ of $\lambda$-small objects of $\mathcal A$ such for every morphism 
$f:S\to\coprod\limits_{i\in I}K_i$ with $S\in\mathcal S$ there are morphisms $g_i:S_i\to K_i$ where $S_i\in\mathcal S$ 
for each $i\in I$ such that $f$ factorizes through 
$$
\coprod\limits_{i\in I}g_i:\coprod\limits_{i\in I}S_i\to\coprod\limits_{i\in I}K_i.
$$
Since these classes are closed under unions, there is the greatest class $\mathcal S$ with this property.  
Its objects are called $\lambda$-\textit{compact}.
}
\end{rem}

\begin{defi}\label{def3.6}
{\em
Let $\ck$ be a category with a zero object $0$. A set $\cg$ of objects is called \textit{weakly generating} 
if $\hom(G,K)=\{0\}$ for each $G\in\cg$ implies that $K=0$.
}
\end{defi}

\begin{rem}\label{re3.7}
{\em
A generating set $\cg$ of objects is clearly weakly generating. Recall that the former concept means that, given two distinct morphisms
$f,g:K_1\to K_2$, there is a morphism $h:G\to K_1$, $G\in\cg$ such that $fh$ and $gh$ are distinct.

M. Hovey proved in \cite{H}, 7.3.1 that the homotopy category of a cofibrantly generated pointed model category has a set
of weak generators.
}
\end{rem}

The following definition is due to A. Neeman.

\begin{defi}\label{def3.8}
{\em
Let $\lambda$ be an infinite cardinal. A category $\ck$ with coproducts and a zero object is called \textit{well}
$\lambda$-\textit{generated} if it has a weakly generating set of $\lambda$-compact objects.

$\ck$ is called \textit{well generated} if it is well $\lambda$-generated for some infinite cardinal $\lambda$.
}
\end{defi}

\begin{theo}\label{th3.9}
Let $\ck$ be a strongly $\lambda$-combinatorial model category. Then $\Ho(\ck)$ is well $\lambda$-generated.
\end{theo}
\begin{proof}
Following \cite{H}, 7.3.1, $\Ho(\ck_\lambda)$ weakly generates $\Ho(\ck)$. Consider a morphism $f:A\to\coprod\limits_{i\in I}K_i$ 
where $A$ is in $\Ho(\ck_\lambda)$. Following \ref{cor3.3}, $E_\lambda f:E_\lambda A\to\coprod\limits_{i\in I}E_\lambda K_i$. 
Since $E_\lambda A$ is $\lambda$-presentable in $\Ind_\lambda\Ho(\ck_\lambda)$ and a coproduct is a $\lambda$-filtered colimit 
of $\lambda$-small subcoproducts, $E_\lambda f$ factorizes through some $\coprod\limits_{j\in J}E_\lambda K_j$ where $J$ has 
the cardinality smaller than $\lambda$. Since $E_\lambda$ is $\Ho(\ck_\lambda)$-full, $A$ is $\lambda$-small.
 
Analogously, the proof of \ref{cor3.3} also yields that objects from $\Ho(\ck_\lambda)$ are $\lambda$-compact. The reason is
that a morphism $f:A\to\coprod\limits_{i\in I}K_i$ with $A\in\Ho(\ck_\lambda)$ is sent by $E_\lambda$ to
the morphism whose codomain is a $\lambda$-filtered colimit of coproducts of objects from $\Ho(\ck_\lambda)$.
\end{proof} 

As a corollary we get the result of A. Neeman \cite{N$_1$} that, for any Grothendieck abelian category $\ck$, 
the derived category $D(\ck)$ is well generated.  

\section{Brown representability}\label{sec5}

\begin{defi}\label{defi4.1}
{\em
A locally $\lambda$-presentable model category $\ck$ will be ca\-lled $\lambda$-\textit{Brown on morphisms} 
provided that the functor 
$$
E_\lambda:\Ho(\ck)\to\Ind_\lambda \Ho(\ck_\lambda)$$ 
is full. $\ck$ will be called $\lambda$-\textit{Brown on objects} provided that $E_\lambda$ is essentially surjective. 
Finally, $\ck$ is $\lambda$-\textit{Brown} if it is $\lambda$-Brown both on objects and on morphisms.
}
\end{defi}

\begin{rem}\label{rem4.2}
{\em
(i) Recall that $E_\lambda$ is essentially surjective if each object in $\Ind_\lambda \Ho(\ck_\lambda)$ is isomorphic
to $E_\lambda K$ for some $K$ in $\Ho(\ck)$. 

(ii) Whenever $\ck$ is strongly $\omega$-combinatorial and $E_\omega$ is full then it is essentially surjective on objects as well. 
In fact, by \ref{cor3.2}, $\Ind_\omega P_\omega$ is full. Since each object of $\Ind_\omega (\ck_\omega)$ can be obtained 
by an iterative taking of colimits of smooth chains (see \cite{AR}) and $P_\omega$ is essentially surjective on objects, 
$\Ind_\omega P_\omega$ is essentially surjective on objects as well. Hence $\ck$ is $\omega$-Brown
on objects. This argument does not work for $\lambda >\omega$ because, in the proof, we need colimits of
chains of cofinality $\omega$. This result corresponds to \cite{Be}, 11.8.

(iii) If $\ck$ is a locally finitely presentable model category such that $\Ho(\ck)$
is a stable homotopy category in the sense of \cite{HPS} then $\ck$ is
$\omega$-Brown in our sense iff $\Ho(\ck)$ is Brown in the sense of \cite{HPS}.

(iv) Let $\ck$ be a strongly $\lambda$-combinatorial model category which is $\lambda$-Brown on morphisms. Consider an object $K$ 
in $\ck$. We can express $K$ as a $\lambda$-filtered colimit of a diagram $D:\cd\to\Ho(\ck_\lambda)$ with a colimit cocone 
$(\delta_d:Dd\to K)_{d\in\cd}$. We get the cone $(P\delta_d:PDd\to PK)_{d\in\cd}$ and, following \ref{cor3.2}, 
$(E_\lambda P\delta_d:E_\lambda PDd\to E_\lambda PK)_{d\in\cd}$ is a colimit cocone. Let $\varphi_d:PDd\to L$ be another cocone. 
There is a unique morphism $t:E_\lambda PK\to E_\lambda L$ such that $tE_\lambda P\delta_d=E_\lambda \varphi_d$ for each $d\in\cd$. 
Since $\ck$ is $\lambda$-Brown on morphisms, we have $t=E_\lambda\overline{t}$ where $\overline{t}:PK\to L$. Since $E_\lambda$ 
is $\Ho(\ck_\lambda)$-faithful, $\overline{t}P\delta_d=\varphi_d$ for each $d\in\cd$. Hence $P\delta_d:PDd\to PK$ is a weak colimit. 
Hence each object of $\Ho(\ck)$ is a weak $\lambda$-filtered colimit of objects from $\Ho(\ck_\lambda)$.

Consider a morphism $f:PK\to PK$ such that $fP\delta_d=P\delta_d$ for each $d$ in $\cd$. Then $E_\lambda f$ is an isomorphism and, 
if $E_\lambda$ reflects isomorphisms, $f$ is an isomorphism as well. This means that each object in $\Ho(\ck)$ is a \textit{minimal 
weak colimit} (in the sense of \cite{HPS}) of objects from $\Ho(\ck_\lambda)$. Minimal colimits are determined uniquely 
up to an isomorphism. Another possible terminology, going back to \cite{Ha}, is a \textit{stable weak colimit}.

(v) Let $\ck$ be a strongly $\lambda$-combinatorial model category which is $\lambda$-Brown. Consider a $\lambda$-filtered diagram 
$D:\cd\to\Ho(\ck_\lambda)$ and let $(\delta_d :E_\lambda Dd \to K)_{d\in\cd}$ be a colimit of $E_\lambda D$ in $\Ind_\lambda\Ho(\ck_\lambda))$. 
Since $\ck$ is $\lambda$-Brown on objects, we can assume that $K=E_\lambda\overline{K}$. Since $E_\lambda$ is $\Ho(\ck_\lambda)$-full
and faithful, there is a cocone $\overline{\delta}_d:Dd\to\overline{K}$ such that $E_\lambda\overline{\delta}_d=\delta_d$ for each
$d$ in $\cd$. By the same argument as in (iv), we get that $\overline{\delta}_d:Dd\to\overline{K}$ is a weak colimit cocone. Hence
$\Ho(\ck)$ has weak $\lambda$-filtered colimits of objects from $\Ho(\ck_\lambda)$.

(vi) $\ck$ being $\lambda$-Brown can be viewed as a weak $\lambda$-accessibility
of $\Ho(\ck)$ because $\Ho(\ck)$ is $\lambda$-accessible with $\Ho(\ck)_\lambda
=\Ho(\ck_\lambda)$ iff  
$$
E_\lambda:\Ho(\ck)\to\Ind_\lambda\Ho(\ck_\lambda)
$$
is an equivalence. This means that $\ck$ is $\lambda$-Brown and $E_\lambda$ is faithful.
But this happens very rarely.
}
\end{rem}

\begin{examples}\label{exa4.3}
{\em
We will show that the homotopy categories
$$
\Ho(\SSet_n)
$$
are finitely accessible
for each $n=1,2,\dots$, i.e., that $E_\omega$ is an equivalence in this case.
Recall that $\SSet_n=\Set^{\BD_n}$ where $\BD_n$ is the category of ordinals
$\{1,2,\dots,n\}$. The model
category structure is the truncation of that on simplicial sets, i.e., cofibrations
are monomorphisms and trivial cofibrations are generated by the horn inclusions
$$
j_m:\Delta^k_m\to\Delta_m\quad\quad 0<k\leq m\leq n.
$$
Here, $\Delta_m=Y_n(m+1)$ where $Y_n:\BD_n\to\SSet_n$ is the Yoneda em\-be\-dding
for $m<n$ and $\Delta_n$ is $Y_{n}(n+1)$ without the $(n+1)$-dimensional simplex
$\{0,1,\dots,n\}$.

For example $\SSet_1=\Set$ and trivial cofibrations are generated by $j_1:1\to 2$. Then weak equivalences are precisely 
mappings between non-empty sets and $\Ho(\SSet_1)$ is the category $\mathbf 2$; all non-empty sets are weakly equivalent. 
$\SSet_2$ is the category of oriented multigraphs with loops. Trivial cofibrations are generated by the embedding $j_1$ of
$$
\bullet 0
$$
to
$$
\xymatrix@1@C=3pc{
0 \bullet \ar[r]& \bullet 1
}
$$
(degenerated loops are not depicted), by the embedding $j_2$ of
$$
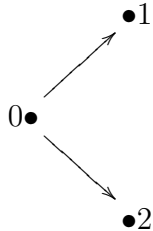
\xymatrix{
& \bullet 1\\
0\bullet \ar[ru] \ar [rd]&\\
& \bullet 2
}
$$
to
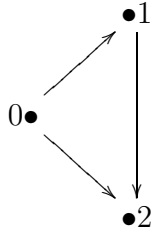
$$
\xymatrix{
& \bullet 1\ar[dd]\\
0\bullet \ar[ru] \ar [rd]&\\
& \bullet 2
}
$$
and their orientation variants. This makes all connected multigraphs weakly equivalent and $\Ho(\SSet_2)$ is equivalent 
to $\Set$; the cardinality of a set corresponds to the number of connected components.

In the case of $\SSet_3$, $1$-connected objects cease to be weakly
equivalent and their contribution to $\Ho(\SSet_3)$ are trees (with a single
root) of height $\,\leq 2$. For example,
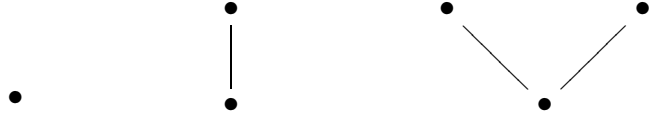
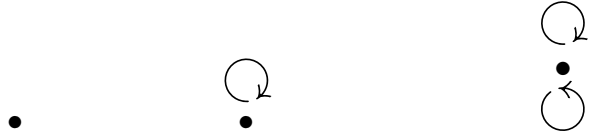
\begin{align*}
&\xymatrix{\\
 \bullet}
\quad\hskip 2cm
\xymatrix{\bullet \ar@{-}[d]\\
\bullet}
\quad\hskip 2cm
\xymatrix{\bullet &&\bullet\\
&\bullet \ar@{-}[ul]\ar@{-}[ur]&
}\\
\intertext{correspond to}
&
\bullet
\quad\hskip 2cm
\stackrel{\mbox{\rotatebox{210}{\Huge
$\circlearrowright$}}}{\bullet}
\quad\hskip 83pt \stackrel{\mbox{\rotatebox{210}{\Huge
$\circlearrowright$}}}{\stackrel{\text{\aa
\char'017}}{\mbox{\rotatebox{10}{\Huge
$\circlearrowleft$}}}}
\end{align*}
\noindent
(degenerated loops are not depicted). Therefore $\Ho(\SSet_3)$ is equivalent
to the category of forests of height $\,\leq 2$. Analogously $\Ho(\SSet_n)$ is
equivalent to the category of forests of height $\,\leq n$. Hence it is
finitely accessible.

Let us add that $\SSet_2$ is a natural model category of oriented multigraphs
with loops (cf. \cite{KR}) and that the symmetric variants
$\Set^{\mathbf F^{op}_n}$, where $\mathbf F_n$ is the category of cardinals
$\{1,\dots,n\}$, are Quillen equivalent to $\SSet_n$ and left-determined by
monomorphisms in the sense of \cite{RT}.
}
\end{examples}

\begin{defi}\label{defi4.4}
{\em
Let\ $\ck$ be a model category. Morphisms $f,g:K\to L$ in $Ho(\ck)$ will be called
$\lambda$-\textit{phantom equivalent} if $E_\lambda f=E_\lambda g$.
}
\end{defi}

This means that $f,g:K\to L$ are $\lambda$-phantom equivalent iff $fh=gh$
for each morphism $h:A\to K$ with $A\in \Ho(\ck_\lambda)$.

\begin{propo}\label{prop4.5}
Let $\ck$ be a strongly $\lambda$-combinatorial model category which is $\lambda$-Brown on morphisms. Then for each object 
$X$ in $\Ho(\ck)$ there exists a weakly initial $\lambda$-phantom equivalent pair $f,g:X\to L$.
\end{propo}

\begin{proof}
We have $X=PK$. Let $(\delta_d:Dd\to K)_{d\in\cd}$ be a canonical $\lambda$-filtered colimit of objects from $\ck_\lambda$.
Following \ref{rem4.2} (iv), $(P\delta_d:PDd\to X)_{d\in\cd}$ is a weak $\lambda$-filtered colimit. Take the induced morphism
$p:\underset{d\in \cd}{\coprod} PDd\to X$ and its weak cokernel pair $f,g$
$$
\xymatrix@C=.2pc@1{
\coprod PDd\  \ar[rrrrr]^{p} &&&&&\ X\ \ar
@<0.6ex>[rrrrrr]^{f}
\ \ar@<-0.6ex>[rrrrrr]_g &&&&&&\  L.
}
$$
Since the starting colimit $(\delta_d:Dd\to K)_{d\in\cd}$ is canonical, $E_\lambda p$ is 
an epimorphism in $\Ind_\lambda \Ho(\ck_\lambda)$. Thus $f$ and $g$ are $\lambda$-phantom equivalent.

Let $f',g':K\to L'$ be a $\lambda$-phantom equivalent. Then $f'p=g'p$ and thus 
the pair $f'$, $g'$ factorizes through $f$, $g$.
Thus $f$, $g$ is a weakly initial $\lambda$-phantom equivalent pair.
\end{proof}

For $\lambda <\mu$ we get a unique functor
$$
F_{\lambda\mu}:\Ind_\mu(\Ho(\ck_\mu))\to\Ind_\lambda(\Ho(\ck_\lambda))
$$
which preserves $\mu$-filtered colimits and whose domain restriction on
$\Ho(\ck_\mu)$ coincides with that of $E_\lambda$.

\begin{propo}\label{prop4.6}
Let $\ck$ be a locally $\lambda$-presentable strongly $\mu$-com\-bi\-na\-to\-rial model where $\lambda <\mu$ are regular cardinals. 
Then $F_{\lambda\mu}E_\mu\cong E_\lambda$.
\end{propo}
\begin{proof}
Following \ref{cor3.2}, we have $E_\mu P\cong\Ind_\mu (P_\mu)$ and thus the functors
$F_{\lambda\mu}E_\mu P\cong F_{\lambda\mu}\Ind_\mu (P_\mu)$ and $E_\lambda P$ have
the isomorphic domain restrictions on $\ck_\mu$. We will show that the functor
$E_\lambda P$ preserves $\mu$-filtered colimits. Since $F_{\lambda\mu}\Ind_\mu (P_\mu)$
has the same property, we will obtain that $F_{\lambda\mu}E_\mu P\cong E_\lambda P$
and thus $F_{\lambda\mu}E_\mu\cong E_\lambda$.

The functor $E_\lambda P$ preserves $\mu$-filtered colimits iff for every object $A$
in $\ck_\lambda$ the functor
$$
hom(PA,P-):\ck\to\Set
$$
preserves $\mu$-filtered colimits. Since $\ck_\lambda\subseteq\ck_\mu$, this follows
from \ref{th3.1}.
\end{proof}

\begin{coro}\label{cor4.7}
Let $\ck$ be a locally $\lambda$-presentable strongly $\mu$-com\-bi\-na\-to\-rial model category where $\lambda <\mu$ 
are regular cardinals. Then $E_\mu$ reflects isomorphisms provided that $E_\lambda$ reflects isomorphisms.
\end{coro}
\begin{proof}
It follows from \ref{prop4.6}.
\end{proof} 

M. Hovey \cite{H} introduced the concept of a pre-triangulated category (distinct from that used in
\cite{N}) and showed that the homotopy category of every pointed model category is pre-triangulated
in his sense. He calls a pointed model category $\ck$ \textit{stable} if $\Ho(\ck)$ is triangulated. In particular, 
$\ck$ is stable provided that $\Ho(\ck)$ is a stable homotopy category in the sense of \cite{HPS}. 

\begin{propo}\label{prop4.8}
Let $\ck$ be a strongly combinatorial stable model category. Then $E_\lambda$ reflects isomorphisms for arbitrarily large regular 
cardinals $\lambda$.
\end{propo}

\begin{proof}
Following \cite{H} 7.3.1, every combinatorial pointed model category $\ck$ has a set $\cg$ of weak generators. 
Let $\Sigma^\ast=\{\Sigma^n Z|Z\in\cg,n\in\mathbf Z\}$. Following \cite{N} 6.2.9, there is a regular 
cardinal $\lambda$ such that $E_\lambda$ reflects isomorphisms. Thus the result follows from \ref{cor4.7}.
\end{proof}

\vskip 0.5cm
\baselineskip 10pt
\noindent
Department of Mathematics\newline
Masaryk University\newline
60000 Brno, Czech Republic\newline
rosicky@math.muni.cz\newline

\end{document}